\documentstyle[amssymb,amsfonts]{amsart}

\def\be#1{\begin{equation} \label{#1}}
\def\bs{\begin{split}}
\def\bi{\begin{itemize}}
\def\ei{\end{itemize}}
\def\es{\end{split}}
\def\ba{\begin{align}}
\def\bas{\begin{align*}}
\def\ea{\end{align}}
\def\eas{\end{align*}}

\def\R{{\mbox{\bf R}}}

\def\eps{\varepsilon}

\def\emph#1{{\it #1}}
\def\textbf#1{{\bf #1}}
\newenvironment{proof}{\noindent {\bf Proof} }{\endprf\par}
\def \endprf{\hfill  {\vrule height6pt width6pt depth0pt}\medskip}

\theoremstyle{plain}
  \newtheorem{theorem}[subsection]{Theorem}
  
  \newtheorem{lemma}[subsection]{Lemma}
  \newtheorem{corollary}[subsection]{Corollary}

\theoremstyle{remark}

\theoremstyle{definition}

\include{psfig}

\begin{document}

\title[A converse extrapolation theorem]{A converse extrapolation theorem for translation-invariant operators }

\author{Terence Tao}
\address{Department of Mathematics, UCLA, Los Angeles, CA 90024}
\email{tao@@math.ucla.edu}

\begin{abstract}
We prove a converse of Yano's extrapolation theorem for translation invariant operators.
\end{abstract}

\maketitle

\section{Introduction}

Let $X$ be a compact symmetric space with compact symmetry group $G$, and $r >  0$, $1 < p_0 < \infty$ be numbers; all constants may depend on $r$ and $p_0$.  If a linear operator $T$ is bounded on $L^p$, $1 < p < p_0$ with an operator norm of $O((p-1)^{-r})$ as $p \to 1$, then it is a classical extrapolation theorem of Yano \cite{yano} that $T$ also maps $L \log^r L(X)$ to $L^1(X)$.

In this paper we show the following converse:

\begin{theorem}\label{main}  Let $G$, $X$, $p_0$, and $r$ be as above.  Suppose $T$ is translation invariant, maps $L \log^r L$ to $L^1$, and is bounded on $L^{p_0}$.  Then $T$ is bounded on $L^p$, $1 < p < p_0$  with an operator norm of $O((p-1)^{-r})$.
\end{theorem}

This theorem is false without the assumption of translation invariance, since $L^p$ is not an interpolation space between $L \log^r L$ and $L^{p_0}$. For a concrete counterexample, take $E$ and $F$ be subsets of $X$ of measure $2^{-N}$ and $N^{rp_0'} 2^{-N}$ respectively, where $N$ is a large number.  Then the operator
$$ Tf = 2^N N^{-r/(p_0-1)} \langle f, \chi_E \rangle \chi_F$$
maps $L \log^r L$ to $L^1$ and bounded on $L^{p_0}$, but the $L^p$ operator norm for $1 < p < p_0$ grows polynomially in $N$.

The translation invariance hypothesis is exploited via the following heuristic principle: if $f$ is a function on $X$ supported on a set of measure $O(1/N)$, then there exists $N$ translates of $f$ which are essentially disjoint.  This idea is used in factorization theory (see e.g. \cite{borg:kakeya}) and also appears in the abstract theory of covering lemmas (e.g. \cite{cordoba:covering}, \cite{guzman:text}).  The point is that the $(L \log^r L, L^1)$ hypothesis yields more information when applied to the sum of the $N$ translates of $f$ than when applied to just $f$ by itself.

The theorem also holds for $p_0 = \infty$, either by a routine modification of the argument, or by assuming an a priori operator bound on $L^2$ (for instance), applying the theorem with $p_0 = 2$, and re-interpolating the result with $L^\infty$ to obtain a better bound on $L^2$.  The theorem also holds of course for $r=0$ by Riesz convexity.

Although our theorem is phrased for compact spaces, it can be extended to non-compact Lie groups if all operator norms are local.  In other words, if $T$ is translation invariant, locally bounded on $L^{p_0}$ and locally maps $L \log^r L$ to $L^1$, then $T$ is also locally bounded on $L^p$, $1 < p < p_0$, with an operator norm of $O((p-1)^{-r})$.  This can be proven either by direct modification of the argument, or by abstract transplantation considerations.

As is well known, the space $L \log^r L$ is an atomic space generated by the atoms $|E|^{-1} \log(1/|E|)^{-r} \chi_E$, where $E$ is an arbitrary measure subset of $X$ with $0 < |E| \ll 1$.  (For completeness, we provide a proof of this fact in an appendix).  As a consequence we have

\begin{corollary}  Let $G$, $X$, $p_0$, $r$ be as above, and let $T$ be a translation invariant operator which is bounded on $L^{p_0}(X)$.  Then a necessary and sufficient condition for $T$ to be bounded on $L^p$, $1 < p < p_0$, with an operator norm of $O(1/(p-1)^r)$, is that
$$ \int |T \chi_E| \lesssim |E| \log(1/|E|)^r$$
for all measurable subsets $E$ of $X$ with $0 < |E| \ll 1$.
\end{corollary}

In a subsequent paper with Jim Wright \cite{tw:multiplier}, we show that certain classes of rough multipliers are bounded from $L \log^r L$ to $L^1$ for various values of $r$, and apply Theorem \ref{main} to deduce sharp bounds for the growth of $L^p$ operator norms.

\section{The main lemma}

We use $A \lesssim B$ to denote the estimate $A \leq CB$ where $C$ is a constant depending on $p_0$, $r$, and the implicit constants in Theorem \ref{main}, and $A \sim B$ to denote the estimates $B \lesssim A \lesssim B$.

Fix $p$; by Riesz convexity we may assume that $p < \frac{1 + p_0}{2}$.  All of our implicit constants shall be independent of $p$.

The main lemma in the argument is

\begin{lemma}\label{key}  Let $E, F$ be subsets of $X$ with $0 < |E| \leq |F|$.  Then we have
\be{chuck} \int_F |Tf| \lesssim |E|^{1/p'} (\frac{1}{p-1} + \log (2+\frac{|F|}{|E|}))^r \|f\|_p
\end{equation}
for all $L^p$ functions $f$ supported on $E$.
\end{lemma}

We remark that without translation invariance, one can only obtain \eqref{chuck} with $\log (2 + \frac{|F|}{|E|})$ replaced by $\log(2 + \frac{1}{|E|})$.

\begin{proof}
Fix $E$, $F$, $f$; we may normalize $\|f\|_p = 1$.  Let $h$ denote the function $h = |\chi_F Tf|$, and define the quantity $A$ by
\be{a-targ} \|h\|_1 = A |E|^{1/p'};
\end{equation}
our task is then to show that 
\be{a-final}
A \lesssim (\frac{1}{p-1} + \log (2 + \frac{|F|}{|E|}))^r.
\end{equation}

Let $N$ be the nearest integer to $\eps/|F|$, where $0 < \eps \ll 1$ is a small constant to be chosen later.  The first step in the argument is to construct group elements $\Omega_0, \ldots, \Omega_N \in G$ such that
\be{cum-b} \bigl \langle \chi_{\bigcup_{j < J} \Omega_j(F)}, h \circ \Omega_J \bigr \rangle \leq \frac{1}{2} A |E|^{1/p'}
\end{equation}
and
\be{cum-a}
 \bigl \langle (\sum_{j < J} |f| \circ \Omega_j)^{p-1}, |f| \circ \Omega_J \bigr \rangle \leq 1
\end{equation}
for all $0 \leq J \leq N$.  

Intuitively, \eqref{cum-b} asserts that the $h \circ \Omega_j$ are essentially disjoint, while \eqref{cum-a} asserts that the $|f| \circ \Omega_j$ are similarly disjoint.  For future reference, we note that \eqref{cum-a} and the $L^p$ normalization of $f$ implies that
\be{recur} \int_X (\sum_{j \leq J} |f| \circ \Omega_j)^p - \int_X (\sum_{j < J} |f| \circ \Omega_j)^p \leq C.
\end{equation}

We now construct the desired group elements.  We may let $\Omega_0$ be arbitrary since \eqref{cum-b}, \eqref{cum-a} are vacuously true for $J=0$.  Now suppose inductively that $\Omega_0, \ldots \Omega_{J-1}$ have already been constructed for some $0 < J \leq N$ such that \eqref{cum-a} (and hence \eqref{recur}) holds for all previous values of $J$.  We will show that
\be{av-b} \int_G \langle \chi_{\bigcup_{j < J} \Omega_j(F)}, h \circ \Omega_J \bigr \rangle\ d\Omega_J \leq \frac{1}{8} A |E|^{1/p}.
\end{equation}
and
\be{av-a} \int_G \bigl \langle (\sum_{j < J} |f| \circ \Omega_j)^{p-1}, |f| \circ \Omega_J \bigr \rangle\ d\Omega_J \leq \frac{1}{4}
\end{equation}
where $d\Omega_J$ is Haar measure on $G$.  By Markov's inequality, this implies that a randomly selected $\Omega_J$ has probability at least 3/4 of obeying \eqref{cum-b} and probability at least 3/4 of obeying \eqref{cum-a}, and so there exists an $\Omega_J$ with the desired properties.

From Fubini's theorem, \eqref{a-targ}, and the identity
$$ \int_G g \circ \Omega(x)\ d\Omega = C \int_X g$$
for all $x \in X$, the left-hand side of \eqref{av-b} evaluates to 
$$C |\bigcup_{j < J} \Omega_j(F)| A |E|^{1/p} \lesssim J |F| A |E|^{1/p}
\lesssim \eps A |E|^{1/p}.$$
Thus \eqref{av-b} holds if $\eps$ is sufficiently small.  The left-hand side of \eqref{av-a} can similarly be evaluated as
$$ C (\int_X (\sum_{j < J} |f| \circ \Omega_j)^{p-1}) (\int_X |f|).$$
From H\"older we have
$$ \int_X |f| \leq |E|^{1/p'} \|f\|_p = |E|^{1/p'},$$
and
$$ \int_X (\sum_{j < J} |f| \circ \Omega_j)^{p-1} \leq (J|E|)^{1/p} 
(\int_X (\sum_{j < J} |f| \circ \Omega_j)^p)^{\frac{p-1}{p}}.$$
On the other hand, from \eqref{recur} and the induction hypothesis we have
$$ \int_X (\sum_{j < J} |f| \circ \Omega_j)^p \lesssim J.$$
Combining all these estimates, we see that 
$$ \hbox{LHS of \eqref{av-a}} \lesssim 
J |E| \lesssim \eps |E|/|F| \lesssim \eps$$
Thus we obtain \eqref{av-a} if $\eps$ is sufficiently small.

Fix $\eps$; all constants may now implicitly depend on $\eps$.  By telescoping \eqref{recur} we have
\be{sa}
\int_X (\sum_{j \leq N} |f| \circ \Omega_j)^p \lesssim N \lesssim |F|^{-1}.
\end{equation}
Let $\epsilon_j = \pm 1$ be an arbitrary assignment of signs.  Then the function
$\sum_{j \leq N} \epsilon_j f \circ \Omega_j$ has a $L^p$ norm of $O(|F|^{-1/p})$ and is supported on a set of measure $O(N|E|) = O(|E|/|F|)$.  We now apply

\begin{lemma}  Let $g$ be a function supported on a set $E \subset E$.  Then
$$ \| g \|_{L \log^r L} \lesssim (\frac{1}{p-1} + \log(2 + \frac{1}{|E|}))^r |E|^{1/p'} \|g\|_p.$$
\end{lemma}

\begin{proof}  We divide into two cases, $|E| \geq 2^{-2r/(p-1)}$
and $|E| \leq 2^{-2r/(p-1)}$.  We normalize
$$ \|g\|_p = (p-1)^r$$
in the first case and
$$ \|g\|_p = (\log \frac{1}{|E|})^{-r} |E|^{-1/p'}$$
in the second; in either case our task reduces to showing that
$$  \int_E |g| \log(2+|g|)^r \lesssim 1.$$
We may restrict ourselves to the set
$$ E' = \{ x \in E: |g(x)| \geq 2 + |E|^{-1} \log^{-r} (2+\frac{1}{|E|}) \},$$
since the contribution outside of $E'$ is clearly acceptable.  In this set $\log(2+|g|)$ may of course be replaced by $\log|g|$.

The function $\frac{\log^r t}{t^{p-1}}$ is increasing for
$1 \leq t < e^{r/(p-1)}$ and decreasing for $t > e^{r/(p-1)}$, with a global maximum of $\frac{(r/e)^r}{(p-1)^r}$. We thus have
$$ \frac{\log(|g|)^r}{|g|^{p-1}} \lesssim \frac{1}{(p-1)^r}$$
in the first case and
$$ \frac{\log(|g|)^r}{|g|^{p-1}} \lesssim |E|^{p-1} \log^{pr} \frac{1}{|E|}$$
if the second case.  In either case the claim follows by multiplying this estimate by $|g|^p$ and integrating, using the $L^p$ normalization of $g$.
\end{proof}

From this lemma we obtain
$$ \| \sum_{j \leq N} \epsilon_j f \circ \Omega_j \|_{L \log^r L} \lesssim
|E|^{1/p'} |F|^{-1} (\frac{1}{p-1} + \log (2+\frac{|F|}{|E|}))^r.$$
Since $T$ is translation invariant and maps $L \log^r L$ to $L^1$, we thus have
$$ \| \sum_{j \leq N} \epsilon_j Tf \circ \Omega_j \|_1 \lesssim
|E|^{1/p'} |F|^{-1} (\frac{1}{p-1} + \log (2+\frac{|F|}{|E|}))^r.$$
Randomizing the signs $\epsilon_j$ and taking expectations using Khinchin's inequality, we obtain
$$ \| (\sum_{j \leq N} |Tf \circ \Omega_j|^2)^{1/2} \|_1 \lesssim
|E|^{1/p'} |F|^{-1} (\frac{1}{p-1} + \log (2+\frac{|F|}{|E|}))^r.$$
In particular, we have
\be{attach} \| (\sum_{j \leq N} (h \circ \Omega_j)^2)^{1/2} \|_1 \lesssim
|E|^{1/p'} |F|^{-1} (\frac{1}{p-1} + \log (2+\frac{|F|}{|E|}))^r.
\end{equation}

If we integrate the trivial pointwise estimate
$$
 (\sum_{j \leq J} (h \circ \Omega_j)^2)^{1/2}
\geq (\sum_{j < J} (h \circ \Omega_j)^2)^{1/2} + h \circ \Omega_J
(1 - \chi_{\bigcup_{j < J} \Omega_j(F)})
$$
using \eqref{a-targ} and \eqref{cum-b}, we obtain
$$
\| (\sum_{j \leq J} (h \circ \Omega_j)^2)^{1/2} \|_1
\geq
\| (\sum_{j < J} (h \circ \Omega_j)^2)^{1/2} \|_1
+ \frac{1}{2} A |E|^{1/p'}.$$
Telescoping this for all $1 \leq J \leq N$, we obtain
$$ \| (\sum_{j \leq N} (h \circ \Omega_j)^2)^{1/2} \|_1 \geq \frac{1}{2} N A |E|^{1/p'} \sim A |E|^{1/p'} |F|^{-1}.$$
Comparing this with \eqref{attach} we obtain \eqref{a-final} as desired.
\end{proof}

\section{Conclusion of the argument}

We are now ready to prove Theorem \ref{main}.  By duality, it suffices to prove the bilinear form estimate
\be{bil}
 |\langle Tf, g \rangle| \lesssim \frac{1}{(p-1)^r}
\end{equation}
for all $f, g$ such that $\|f\|_p = 1$, $\|g\|_{p'} = 1$.

Fix $f$, $g$; we may assume that $f$, $g$ are non-negative.  Let $f^*: \R^+ \to \R^+$ be the non-increasing left-continuous re-arrangement of $f$, so that $\|f^*\|_p = 1$ and
\be{flip} | \{ x: f(x) > f^*(\alpha) \} | \leq \alpha.
\end{equation}
Similarly define $g^*$.

For any integers $q \geq 1$ and $k < C$, define $f_{k,q}$ to be the restriction of $f$ to the set $\{x: f^*(2^{qk+q}) < f(x) \leq f^*(2^{qk})\}$.  Since $X$ has finite measure, we thus have
$$ f = \sum_{k} f_{k,q}.$$
Similarly define $g_{k,q}$.

Usually one takes $q = 1$, but because of our desire for sharp bounds as $p \to 1$ it shall be more appropriate to choose $q$ so that $q \sim 1/(p-1) \sim p'$.  

By the triangle inequality, \eqref{bil} will now follow from the estimates
\be{bil2}
 \sum_{k, l: k \geq l+1} |\langle Tf_{k,q}, g_{l,q} \rangle| \lesssim 1
\end{equation}
and
\be{bil3}
 \sum_{k, l: k \leq l} |\langle Tf_{k,q}, g_{l,q} \rangle| \lesssim q^r.
\end{equation}

Let us first prove \eqref{bil2}.  By splitting
$$ f_{k,q} = \sum_{k' = qk}^{qk+q-1} f_{k',1}, \quad g_{l,q} = \sum_{l' = ql}^{ql+q-1} g_{l',1}$$
we see that the left-hand side of \eqref{bil2} is majorized by
$$ \sum_{k',l': k' > l'} |\langle Tf_{k',1}, g_{l',1} \rangle|.$$
Since $T$ is bounded on $L^{p_0}$, we may use H\"older's inequality to majorize this by
$$ \sum_{k',l': k' > l'} \|f_{k',1}\|_{p_0} \|g_{l',1}\|_{p'_0}.$$
From \eqref{flip} and the definition of $f_{k,q}$ we have
$$ \|f_{k',1}\|_{p_0} \lesssim 2^{k'/p_0} f^*(2^{k'})$$
and similarly
$$ \|g_{l',1}\|_{p_0} \lesssim 2^{l'/p'_0} g^*(2^{l'}).$$
Thus the left-hand side of \eqref{bil2} is majorized by
$$ \sum_{k',l': k' > l'} (2^{k'/p} f^*(2^{k'})) (2^{l'/p'} g^*(2^{l'})) 2^{-(k'-l')(\frac{1}{p} - \frac{1}{p_0})}.$$
The estimate \eqref{bil2} then follows from Young's inequality for bilinear forms.  Indeed, the first expression in parentheses has an $l^p$ norm comparable to $\|f^*\|_p = 1$, the second expression in parentheses has an $l^{p'}$ norm comparable to $\|g^*\|_{p'} = 1$, and convolution kernel is summable with $l^1$ norm of $O(1)$ since we are assuming $p < (1+p_0)/2$.  

It remains to prove \eqref{bil3}.  From Lemma \ref{key} and \eqref{flip} we have
$$
\int_{g > g^*(2^{ql+q})} |Tf_{k,q}| \lesssim (2^{qk+q})^{1/p'} 
(\frac{1}{p-1} + \log (2+\frac{2^{ql+q}}{2^{qk+q}}))^r \|f_{k,q}\|_p.$$
From the definition of $q$ and the assumptions on $k$, $l$, this simplifies to
$$
\int_{g > g^*(2^{ql+q})} |Tf_{k,q}| \lesssim 2^{qk/p'}
q^r (1 + l-k)^r \|f_{k,q}\|_p.$$
From H\"older's inequality we thus have
$$ |\langle Tf_{k,q}, g_{k,q} \rangle| \lesssim g^*(2^{ql}) 2^{qk/p'} q^r (1 + l-k)^r \|f_{k,q} \|_p.$$
Thus the left-hand side of \eqref{bil3} is majorized by
$$ q^r \sum_{k,l: k \leq l} \|f_{k,q}\|_p (2^{ql/p'} g^*(2^{ql})) 
2^{-q(l-k)/p'} (1 + l-k)^r.$$
The claim then follows again from Young's inequality and the choice of $q$, since the sequence $\|f_{k,q}\|_p$ is in $l^p$, the sequence $2^{-ql/p'} g^*(2^{ql})$ has an $l^{p'}$ norm comparable to $\|g^*\|_{p'} = 1$, and the convolution kernel is integrable uniformly in $p$.
\endprf

\section{Appendix: atomic decomposition of Orlicz spaces}

In this section we show that every $L \log^r L(X)$ function $f$ can be decomposed into a convex linear combination of atoms 
$|E|^{-1} \log(1/|E|)^r \chi_E$ with $0 < |E| \ll 1$.

We first observe that any function $f$ supported on a set of measure $2^{-k}$ and having a sup norm of $k^{-r} 2^{-k}$ can easily be decomposed in this manner, since bounded functions can be written as convex linear combinations of characteristic functions.

Now let $f$ be a general $L \log^r L(X)$ function; we may normalize so that
$$ \int_X |f| \log^r(2 + |f|) = 1.$$
We may also assume without loss of generality that $f$ is non-negative and is supported on a set of measure $\ll 1$.

Let $f^*$ and $f_{k,q}$ be as before.  For each integer $k < -C$, we define
$$ c_k = |k|^r 2^k f^*(2^k)$$
and
$$ a_k(x) = f_{k,1} / c_k.$$
Clearly $f = \sum_{k < -C} c_k a_k$.  From \eqref{flip} and the previous discussion, the $a_k$ are convex linear combinations of atoms uniformly in $k$, so it suffices to show that the $c_k$ are summable, i.e. that
$$ \sum_{k < -C} |k|^r 2^k f^*(2^k) \lesssim 1.$$

Since $f$ is non-decreasing, we may bound this expression by
$$ C \int_{0 \leq t \ll 1} f^*(t) \log(1/t)^{-r}\ dt.$$
The portion of the integral where $f^*(t) \leq t^{-1/2}$ is clearly acceptable, so we may assume that $f^*(t) \geq t^{-1/2}$.  But then we may estimate the above by
$$ C \int_{0 \leq t \ll 1} f^*(t) \log(2 + f^*(t))^{-r}\ dt =
C \int f \log^r(2 + f) = C$$
as desired.


\begin{thebibliography}{10}

\bibitem{borg:kakeya}
J. Bourgain, \emph{Besicovitch-type maximal operators and
applications to Fourier analysis}, Geom. and Funct. Anal. \textbf{22}
(1991), 147--187.

\bibitem{cordoba:covering}
A. C\'ordoba, \emph{Maximal functions, covering lemmas and Fourier multipliers}, Harmonic analysis in Euclidean spaces (Proc. Sympos. Pure Math., Williams Coll., Williamstown, Mass., 1978), Part 1, pp. 29--50, Proc. Sympos. Pure Math., XXXV, Part, Amer. Math. Soc., Providence, R.I., 1979.

\bibitem{guzman:text}
M. de Guzm\'an, \emph{Real variable methods in Fourier analysis}, North-Holland Mathematics Studies, 46. Notas de Matemática [Mathematical Notes], 75. North-Holland Publishing Co., Amsterdam-New York, 1981.

\bibitem{tw:multiplier} T. Tao, J. Wright, \emph{Endpoint multiplier theorems of Marcinkiewicz type}, submitted.

\bibitem{yano} S. Yano, \emph{Notes on Fourier analysis. XXIX. An extrapolation theorem. }
J. Math. Soc. Japan \textbf{3}, (1951). 296--305. 


\end{thebibliography}
\end{document}